# An Algebraic Identity Leading to Wilson's Theorem

Sebastián Martín Ruiz
Avda. de Regla 43 Chipiona 11550 Spain
smruiz@telefonica.net

In most text books on number theory Wilson's theorem is proved by applying Lagrange's theorem concerning polynomial congruences [1,2,3,4]. Hardy and Wright also give a proof using quadratic residues [3]. In this note Wilson's theorem is derived as a corollary to an algebraic identity.

*Theorem 1:*

For all integers $n \geq 0$ and for all real numbers $x$

$$\sum_{i=0}^{n}(-1)^{i}\binom{n}{i}(x-i)^{n} = n!$$

*Proof:*

We proceed by induction. Let

$$f_n(x) = \sum_{i=0}^{n}(-1)^{i}\binom{n}{i}(x-i)^{n}$$

It is easy to show that $f_0(x) = 1 = 0!$
Assume $f_k(x) = k! \ \forall \ x \in \mathbb{R}$
We consider

$$f'_{k+1}(x) = \frac{d}{dx}\left[\sum_{i=0}^{k+1}(-1)^{i}\binom{k+1}{i}(x-i)^{k+1}\right] =$$

$$= (k+1)\sum_{i=0}^{k+1}(-1)^{i}\frac{(k+1)!}{i!(k+1-i)!}(x-i)^{k}$$

Splitting off the $i = k+1$ term and using

$$\frac{(k+1)!}{i!(k+1-i)!} = \frac{k!}{i!(k-i)!} \cdot \frac{k+1}{k+1-i} = \binom{k}{i}\left(1 + \frac{i}{k+1-i}\right)$$

we get:

$$f'_{k+1}(x) = (k+1)\left[\sum_{i=0}^{k}(-1)^i \binom{k}{i}(x-i)^k + (-1)^{k+1}[x-(k+1)]^k\right] +$$

$$+ (k+1)\left[\sum_{k=1}^{k}(-1)^i \frac{k!(x-i)^k}{(i-1)!(k+1-i)!}\right]$$

$$= (k+1)\left[f_k(x) + \sum_{j=0}^{k}(-1)^{j+1}\binom{k}{j}(x-1-j)^k\right] \text{ where } j = i-1$$

$$= (k+1)[f_k(x) - f_k(x-1)] = (k+1)[k!-k!] = 0.$$

Thus $f_{k+1}(x)$ is constant and in particular

$$f_{k+1}(x) = f_{k+1}(k+1) =$$

$$= \sum_{i=0}^{k}(-1)^i \frac{(k+1)!}{i!(k-i)!}(k+1-i)^k$$

Since the $i = k+1$ term is zero

$$f_{k+1}(x) = (k+1)\sum_{i=0}^{k}(-1)^i \binom{k}{i}(k+1-i)^k =$$

$$= (k+1)f_k(k+1) = (k+1)k! = (k+1)!$$

so $\quad f_k(x) = k! \Rightarrow f_{k+1}(x) = (k+1)!$

Therefore $\quad f_n(x) = n! \; \forall x \in \mathrm{R} \; \forall n \in \mathrm{N}$

*Corollary 1: (Wilson's theorem)*
  For any prime number $p$, we have: $(p-1)! \equiv p-1 \pmod{p}$.

*Proof:*

Let $n = p-1$ where $p$ is prime. We use the formula for $x = 0$.

$$\sum_{i=0}^{p-1}(-1)^i \binom{p-1}{i}(-i)^{p-1} = (p-1)! \quad (1)$$

But for $p > 1 \geq 1$:

$$0 \equiv \binom{p}{i} = \binom{p-1}{i-1} + \binom{p-1}{i} \pmod{p}$$

Thus

$$\left.\begin{array}{rcl}\binom{p-1}{i} &\equiv& -\binom{p-1}{i-1} \pmod{p} \\ \binom{p-1}{0} &\equiv& 1 \pmod{p}\end{array}\right\} \Rightarrow \binom{p-1}{i} \equiv (-1)^i \pmod{p}$$

Using this result in (1) we obtain

$$\sum_{i=0}^{p-1}(-1)^i(-1)^i(-i)^{p-1} \equiv (p-1)! \pmod{p}$$

If $p > 2$, $p$ is odd and therefore $p-1$ is even. Thus we have the relation $(-i)^{p-1} = i^{p-1}$ which allows us to obtain

$$\sum_{i=0}^{p-1} i^{p-1} \equiv (p-1)! \pmod{p} \quad (2)$$

For on the other hand since $p$ is not a factor of $i$, and using Fermat's theorem, we have

$$i^{p-1} \equiv 1 \pmod{p} \quad (3)$$

Combining (2) and (3), we can conclude:

$$\sum_{i=1}^{p-1} 1 \equiv (p-1)! \pmod{p}$$

and this last relation can be written in the form

$$(p-1)! \equiv p-1 \pmod{p}$$

*Corollary 2:*

For all integers $n \geq 0$ and for all real numbers $x$

$$\sum_{i=0}^{n}(-1)^i \binom{n}{i}(x-i)^{n-j} = 0 \quad 1 \leq j \leq n$$

*Proof:*

We consider the algebraic identity of theorem 1

$$f_n(x) = \sum_{i=0}^{n}(-1)^i \binom{n}{i}(x-i)^n$$

Differentiating j times

$$f_n^{(j)}(x) = n(n-1)\cdots[n-(j-1)]\sum_{i=0}^{n}(-1)^i \binom{n}{i}(x-i)^{n-j} = 0$$

Then,

$$\sum_{i=0}^{n}(-1)^i \binom{n}{i}(x-i)^{n-j} = 0 \quad 1 \leq j \leq n$$

*Author:*

Sebastián Martín Ruiz

http://personal.telefonica.terra.es/web/smruiz/

email: smruiz(AT)telefonica.net

Avda. de Regla 43, Chipiona 11550 Spain